\documentclass[11pt]{amsart}

\usepackage{amsfonts,amsmath,amscd,amssymb,amsbsy,amsthm,amstext,amsopn}
\usepackage{fullpage,mathrsfs,subfigure}
\usepackage{pstricks,pst-node,pst-text,pst-tree}
\usepackage{tikz}
\usepackage{stmaryrd}  
\usepackage{extarrows}
\usepackage[hmargin=3cm,vmargin=3.5cm]{geometry}
\usepackage{graphicx}
\usepackage{url}
\usepackage{mathdots}
\usepackage[all,cmtip]{xy}
\usepackage{color}
\usepackage{url}
\usepackage[capitalize]{cleveref}

\newtheorem{theorem}{Theorem}[section]

\newtheorem{acknowledgements}[theorem]{Acknowledgments}

\theoremstyle{definition}

\theoremstyle{remark}

\newcommand{\kk}{\ensuremath{\Bbbk}} 
\newcommand{\CC}{\ensuremath{\mathbb{C}}}
 
\newcommand{\PP}{\ensuremath{\mathbb{P}}}

\newcommand{\ZZ}{\ensuremath{\mathbb{Z}}}

\newcommand{\tot}{\operatorname{tot}}

 \newcommand{\Cl}{\operatorname{Cl}}

\newcommand{\End}{\operatorname{End}}
\newcommand{\Ext}{\operatorname{Ext}}

\newcommand{\Hom}{\operatorname{Hom}}

\newcommand{\Pic}{\operatorname{Pic}}

\newcommand{\Dcal}{\mathcal{D}}

\newcommand{\Ocal}{\mathcal{O}}

\newcommand{\Lcal}{\mathcal{L}}

\newcommand{\Dx}{\mathcal{D}^b(X)}

\renewcommand\mod{\operatorname{mod \textendash \!}}

\begin{document}
\title{QuiversToricVarieties: a package to construct quivers of sections on complete toric varieties}
\author{Nathan Prabhu-Naik}
\date{9th October 2014} 

\begin{abstract}
Given a collection of line bundles on a complete toric variety, the \emph{Macaulay2} package \emph{QuiversToricVarieties} contains functions to construct its quiver of sections and check whether the collection is strong exceptional. It contains a database of full strong exceptional collections of line bundles for smooth Fano toric varieties of dimension less than or equal to 4.
\end{abstract}
\maketitle
\section{Introduction}
\noindent For a collection of non-isomorphic line bundles $\Lcal = \lbrace \Lcal_0 := \Ocal_X, \Lcal_1, \ldots, \Lcal_r \rbrace$ on a complete normal toric variety $X$, the endomorphism algebra $\End(\bigoplus_i \Lcal_i)$ can be described as the quotient of the path algebra of its quiver of sections by an ideal of relations determined by labels on the arrows in the quiver \cite{CrSm}. The vertices of the quiver correspond to the line bundles and there is a natural order on the vertices defined by $i < j$ if $\Hom(L_j, L_i) = 0$. For $i < j$, the number of arrows from $i$ to $j$ is equal to the dimension of the cokernel of the map 
\begin{equation}
\bigoplus_{i<k<j} \Hom(L_i , L_k) \otimes \Hom(L_k , L_j) \longrightarrow \Hom (L_i , L_j),
\end{equation}
and we label each arrow by the toric divisors corresponding to the sections in a basis for the cokernel.
Using the given order on $\Lcal$, the collection is \emph{strong exceptional} if 
\begin{equation}
\Ext^i (L_j, L_k) = 0, \forall j, \ k, \text{ and } i \neq 0.
\end{equation}

\noindent Let $\Dx$ denote the bounded derived category of coherent sheaves on $X$. The collection $\Lcal$ is \emph{full}, or \emph{generates} $\Dx$, if the smallest triangulated full subcategory of $\Dx$ containing $\Lcal$ is $\Dx$ itself. A tilting bundle $T$ on $X$ is a vector bundle such that $T$ generates $\Dx$ and $\Ext^i(T,T) = 0$ for $i>0$; given a full strong exceptional collection of line bundles $\Lcal$ on $X$, the direct sum $\bigoplus_{L_i \in \Lcal} L_i$ is a tilting bundle. The following theorem by Baer and Bondal allows us to understand $\Dx$ in terms of the module category of a finite dimensional algebra. 

\begin{theorem}\cite{Baer,Bond}
Let $T$ be a tilting bundle on $X$, $A = \End(T)$ and $\Dcal^b(\mod A)$ be the bounded derived category of finitely generated right $A$-modules. Then
\begin{equation}
 \mathbf{R}\Hom (T, - ) \colon \Dx \rightarrow \Dcal^b(\mod A)
\end{equation} 
is an equivalence of triangulated categories.
\end{theorem}

\noindent A complete normal toric variety induces a short exact sequence of abelian groups
\begin{equation}
 \begin{CD}\label{ses}   
    0@>>> M @>>> \ZZ^{\Sigma(1)} @>\deg >> \Cl(X)@>>> 0,
\end{CD}
\end{equation}
where $M$ is the character lattice of the dense torus in $X$, $\Sigma(1)$ is the set of rays in the fan $\Sigma$ of $X$, and the map $\deg$ sends a toric divisor $D \in \ZZ^{\Sigma(1)}$ to the rank $1$ reflexive sheaf $\Ocal_{X} (D)$ in the class group $\Cl (X)$ (see for example \cite{Fult}). Showing that $\Lcal$ is strong exceptional in this situation is equivalent to checking that $H^i(X, L^{-1}_j \otimes L_k) =  0$ for $i > 0,\  0 \leq j,k \leq r$. Using a theorem of Eisenbud, Musta{\c{t}}{\u{a}} and Stillman \cite{EiMuSt}, we can determine if the cohomology of $\Ocal_X(D)$ vanishes by considering when $\Ocal_X(D)$ avoids certain affine cones constructed in $\Cl(X)$, which we call \emph{non-vanishing cohomology cones}. The purpose of the package \emph{QuiversToricVarieties} for \emph{Macaulay2} \cite{M2} is to construct the quiver of sections for a collection of line bundles on a complete toric variety and check if the collection is strong exceptional. We note that there does exist computer programs that check if a collection of line bundles on a toric variety is strong exceptional; see for example Perling's \emph{TiltingSheaves} \cite{Perl}.  

Restricting our attention to smooth toric Fano varieties, toric divisorial contractions give the collection of $n$-dimensional toric Fano varieties a poset structure, described for $n=3$ by \cite{Oda} and $n=4$ by \cite{Sato} (see also \cite[Remark 2.4]{Prna}). The contractions induce lattice maps between the short exact sequences (\ref{ses}) determined by the varieties and these lattice maps are an essential ingredient in the proof that each smooth toric Fano variety of dimension $\leq 4$ has a full strong exceptional collection of line bundles \cite[Theorem 6.4]{Prna}. The package \emph{QuiversToricVarieties} contains a database of these lattice maps and of full strong exceptional collections of line bundles on all smooth toric Fano varieties of dimension $\leq 4$. 

In the case when $X$ is a smooth toric Fano variety, let $Y = \tot(\omega_X)$ be the total space of the canonical bundle on $X$. The package \emph{QuiversToricVarieties} contains methods to check if the pullback of a full strong exceptional collection of line bundles on $X$ along the morphism $Y \rightarrow X$ is a tilting bundle on $Y$.

\emph{QuiversToricVarieties} depends on the package \emph{NormalToricVarieties} for the construction of toric varieties and for the database of smooth toric Fano varieties. All varieties are defined over $\kk = \CC$.   

\section{Overview of the Package}

\noindent Let $X$ be a complete normal toric variety constructed in \emph{NormalToricVarieties} with a torsion-free class group. The class group lattice of $X$ has a basis determined by \verb+fromWDivToCl+ and the function \verb+fromPicToCl+ can be used to determine which vectors in the lattice correspond to line bundles. The input for the method \verb+quiver+ is a complete normal toric variety with a torsion-free class group, together with a list of vectors $v_i$ in the class group lattice that correspond to the line bundles $L_i$. The vectors are ordered by \verb+quiver+ and the basis of $\Hom (L_i,L_j)$ is calculated by determining the basis of the multidegree $v_j-v_i$ over the Cox ring of the variety. From this basis, the irreducible maps are chosen and listed as arrows, with the corresponding toric divisors as labels. If some of the vectors do not correspond to line bundles then a quiver is still constructed but the resulting path algebra modulo relations may not be isomorphic to $\End(\bigoplus_{i \in Q_0} E_i)$, where $E_i$ are the rank $1$ reflexive sheaves corresponding to $v_i$. Alternatively, we can produce a quiver by explicitly listing the vertices, the arrows with labels and the variety. The methods \verb+source+, \verb+target+, \verb+label+ and \verb+index+ return the specific details of an arrow in the quiver, a list of which can be accessed by inputting \verb+Q_1+.

Besides the method \verb+quiver+, the method \verb+doHigherSelfExtsVanish+ forms the core of the package. The primary input is a quiver of sections. The method creates the non-vanishing cohomology cones in the class group lattice for $X$ and determines if the vectors $v_i - v_j$ avoid these cones. The cones are determined by certain subsets $I$ of the rays of the fan $\Sigma$ for $X$; if the complement of the supporting cones for $I$ in $\Sigma$ has non-trivial reduced homology, then $I$ is called a \emph{forbidden set} and it determines a cone in $\ZZ^{\Sigma(1)}$. The forbidden sets can be calculated using the function \verb+forbiddenSets+, and the image of a cone determined by a forbidden set under the map \verb+fromWDivToCl X+ is a non-vanishing cohomology cone in $\Cl(X)$.          

A database in \emph{NormalToricVarieties} contains the smooth toric Fano varieties up to dimension $6$ and can be accessed using \verb+smoothFanoToricVariety+. The divisorial contractions between the smooth toric Fano varieties up to dimension $4$ are listed under the \verb+contractionList+ command, and the induced maps between their respective short exact sequences (\ref{ses}) are recalled from a database in \emph{QuiversToricVarieties} using the \verb+tCharacterMap+, \verb+tDivisorMap+ and the \verb+picardMap+ commands. Note that as each variety considered is smooth, its class group is isomorphic to its Picard group.  

The database containing full strong exceptional collections of line bundles for smooth Fano toric varieties in dimension $\leq 4$ can be accessed using \verb+fullStrExcColl+. The collections for the surfaces were calculated by King \cite{King}, the threefolds by Costa--Mir\'{o}-Roig \cite{CoMR1}, Bernardi--Tirabassi \cite{BeTi} and Uehara \cite{Ueha} and the fourfolds by Prabhu-Naik \cite{Prna}. 

\section{An Example}
\noindent We illustrate the main methods in \emph{QuiversToricVarieties} using the blowup of $\PP^2$ at three points, the birationally-maximal smooth toric Fano surface. It is contained in the toric Fano database in \emph{NormalToricVarieties}, which is loaded by the \emph{QuiversToricVarieties} package. 
\begin{verbatim}
i1 : loadPackage "QuiversToricVarieties";
i2 : X = smoothFanoToricVariety(2,4);
\end{verbatim}

\noindent A full strong exceptional collection $\Lcal$, first considered by King \cite{King}, can be recalled from the database and its quiver of sections can be created.
\begin{verbatim}
i3 : L = fullStrExcColl(2,4);
o3 = {{0,0,0,0},{0,0,1,1},{0,1,0,0},{0,1,1,0},{1,0,0,0},{1,0,0,1}}
i4 : Q = quiver(L,X);
\end{verbatim}

\noindent We can view the details of the quiver, either by displaying the arrows at each vertex, or by listing all of the arrows and considering their source, target and label. 
\begin{verbatim}
i5 : Q#0
o5 = HashTable{1 => {x_0x_1 , x_3x_4 }     }
               2 => {x_1x_2 , x_4x_5 }
               3 => {x_2x_3 , x_0x_5 }
               degree => {0, 0, 0, 0}
               
i6 : first Q_1
o6 = arrow_1

i7 : source oo, target oo, label oo
o7 = (0, 1, x_0x_1 )                           
\end{verbatim}

\noindent The forbidden sets of rays can be computed and the collection of line bundles can be checked to be strong exceptional. The method \verb+doHigherSelfExtsVanish+ creates a copy of the non-vanishing cohomology cones in the cache table for $X$, where the cones are given by a vector and a matrix $\{w,M\}$ encoding the supporting closed half spaces of the cone, in which case the lattice points of the cone are $\{ v \in \Cl (X) \mid  M v \leq w \}$. The non-vanishing cone for $H^2$ is displayed below.
\begin{verbatim}
i8 : peek forbiddenSets X
o8 = MutableHashTable{1 => {{0,2},{0,3},{1,3},{0,1,3},{0,2,3},{0,4},{1,4},...}
                      2 => {{0,1,2,3,4,5}}

i9 : doHigherSelfExtsVanish Q
o9 = true

i10 : X.cache.cones#2
o10 = {{| -1 |, | 1 1 1 0 1 |}}
        | -1 |  | 1 0 1 1 1 |
        | -1 |  | 0 1 1 0 0 |
        | -1 |  | 0 0 0 1 1 |                      
\end{verbatim}
Consider the chain of divisorial contractions $X =: X_4 \rightarrow X_3 \rightarrow X_2 \rightarrow X_0$ from $X$ to the toric Fano surfaces numbered $3$, $2$ and $0$ in the database. The contractions induces lattice maps $\Pic(X_4) \rightarrow \Pic(X_3) \rightarrow \Pic(X_2) \rightarrow \Pic(X_0)$ and the method \verb+doHigherExtsVanish+ can check if the non-isomorphic line bundles in the image of $\Lcal$ under these lattice maps are strong exceptional for each contraction.
\begin{verbatim}
i11 : doHigherSelfExtsVanish(Q,{4,3,2,0})
o11 = true
\end{verbatim}
Now consider the morphism $\pi \colon \tot(\omega_X) \rightarrow X$. The pullback $\pi^* (\bigoplus_{L_i \in \Lcal} L_i)$ is a tilting bundle on $Y = \tot (\omega_X)$ if 
\[ H^k(X,L_i \otimes L_j^{-1} \otimes \omega_X^{-m}) = 0\]
for all $k>0$, $m \geq 0$ and $L_i,L_j \in \Lcal$ (see for example \cite[Theorem 6.7]{Prna}). As $\omega^{-1}_X$ is ample, there exists a non-negative integer $n$ such that $L_i \otimes L_j^{-1} \otimes \omega_X^{-m}$ is nef for $0 \leq i,j \leq r$ and $m \geq n$, and hence $H^k(X,L_i \otimes L_j^{-1} \otimes \omega_X^{-m}) = 0$ for all $k>0$ by Demazure vanishing. The method \verb+bundlesNefCheck+ checks for a given integer $n$ whether $L_i \otimes L_j^{-1} \otimes \omega_X^{-n}$ is nef for all $L_i,L_j \in \Lcal$.
\begin{verbatim}
i12 : n=2;
i13 : bundlesNefCheck(Q,n)
o13 = true
\end{verbatim}  If an integer $p$ is included as an additional input in \verb+doHigherSelfExtsVanish+, then the method checks that for all $0 \leq m \leq p$, whether the line bundles $L_i \otimes L_j^{-1} \otimes \omega^{-m}$ avoid the non-vanishing cohomology cones. Note that for our example, the computation above implies that it is enough to use the integer $n-1$. 
\begin{verbatim}
i14 : doHigherSelfExtsVanish(Q,n-1)
o14 = true
\end{verbatim}
For $t \in \{4,3,2,0\}$, let $\{ L_{i,t} \}_{i \in I_t}$ denote the list of non-isomorphic line bundles in the image of $\Lcal$ under the map given by \verb+picardMap+ from $\Pic(X) \rightarrow \Pic(X_t)$, where $I_t$ is an index set. By including the list of divisorial contractions as an input in \verb+doHigherSelfExtsVanish+, we can check that
\[ H^k(X_t,L_{i,t} \otimes (L_{j,t})^{-1} \otimes \omega_{X_t}^{-m}) = 0\]
for $k >0$, $0 \leq m \leq n-1$, $t \in \{4,3,2,0\}$ and all $i,j \in I_t$.
\begin{verbatim}
i15 : doHigherSelfExtsVanish(Q,{4,3,2,0},n-1)
o15 = true
\end{verbatim}
For all $n$-dimensional smooth toric Fano varieties, $1 \leq n \leq 3$, and 88 of the 124 smooth toric Fano fourfolds, the database contains a chain complex of modules over the Cox ring for the variety. The chain complexes are used in \cite{Prna} to show that the collections of line bundles in the database for these varieties are full.
\begin{verbatim}
i16 : C = resOfDiag(2,4);
i17 : SS = ring C;
i18 : C
        6       12       6
o18 = SS  <-- SS   <-- SS
\end{verbatim}
 
\begin{acknowledgements}
The author is extremely grateful to Gregory G. Smith, who developed a substantial part of this package. He wishes to thank Alastair Craw for supervising the project, and Thomas Hawes for his helpful advice. 
\end{acknowledgements}
\bibliography{QuiversToricVarieties}

\begin{thebibliography}{CMR04}

\bibitem[Bae88]{Baer}
Dagmar Baer.
\newblock Tilting sheaves in representation theory of algebras.
\newblock {\em Manuscripta Math.}, 60(3):323--347, 1988.

\bibitem[Bon90]{Bond}
A.~I. Bondal.
\newblock Helices, representations of quivers and {K}oszul algebras.
\newblock In {\em Helices and vector bundles}, volume 148 of {\em London Math.
  Soc. Lecture Note Ser.}, pages 75--95. Cambridge Univ. Press, Cambridge,
  1990.

\bibitem[BT09]{BeTi}
Alessandro Bernardi and Sofia Tirabassi.
\newblock Derived categories of toric {F}ano 3-folds via the {F}robenius
  morphism.
\newblock {\em Matematiche (Catania)}, 64(2):117--154, 2009.

\bibitem[CMR04]{CoMR1}
L.~Costa and R.~M. Mir{\'o}-Roig.
\newblock Tilting sheaves on toric varieties.
\newblock {\em Math. Z.}, 248(4):849--865, 2004.

\bibitem[CS08]{CrSm}
Alastair Craw and Gregory~G. Smith.
\newblock Projective toric varieties as fine moduli spaces of quiver
  representations.
\newblock {\em Amer. J. Math.}, 130(6):1509--1534, 2008.

\bibitem[EMS00]{EiMuSt}
David Eisenbud, Mircea Musta{\c{t}}{\v{a}}, and Mike Stillman.
\newblock Cohomology on toric varieties and local cohomology with monomial
  supports.
\newblock {\em J. Symbolic Comput.}, 29(4-5):583--600, 2000.
\newblock Symbolic computation in algebra, analysis, and geometry (Berkeley,
  CA, 1998).

\bibitem[Ful93]{Fult}
William Fulton.
\newblock {\em Introduction to toric varieties}, volume 131 of {\em Annals of
  Mathematics Studies}.
\newblock Princeton University Press, Princeton, NJ, 1993.
\newblock The William H. Roever Lectures in Geometry.

\bibitem[GS]{M2}
Daniel~R. Grayson and Michael~E. Stillman.
\newblock Macaulay2, a software system for research in algebraic geometry.
\newblock Available at \url{http://www.math.uiuc.edu/Macaulay2/}.

\bibitem[Kin97]{King}
Alastair King.
\newblock Tilting bundles on some rational surfaces.
\newblock Available as a preprint at
  \url{http://www.maths.bath.ac.uk/~masadk/papers}, 1997.

\bibitem[Oda88]{Oda}
Tadao Oda.
\newblock {\em Convex bodies and algebraic geometry}, volume~15 of {\em
  Ergebnisse der Mathematik und ihrer Grenzgebiete (3) [Results in Mathematics
  and Related Areas (3)]}.
\newblock Springer-Verlag, Berlin, 1988.
\newblock An introduction to the theory of toric varieties, Translated from the
  Japanese.

\bibitem[Per]{Perl}
Markus Perling.
\newblock Tiltingsheaves, a program for computing exceptional sequences on
  toric varieties.
\newblock Available at
  \url{http://www.mathematik.uni-kl.de/~perling/ts/ts.html}.

\bibitem[PN15]{Prna}
Nathan Prabhu-Naik.
\newblock Tilting bundles on toric fano fourfolds.
\newblock Available as a preprint at \url{http://arxiv.org/abs/1501.05871},
  2015.

\bibitem[Sat00]{Sato}
Hiroshi Sato.
\newblock Toward the classification of higher-dimensional toric {F}ano
  varieties.
\newblock {\em Tohoku Math. J. (2)}, 52(3):383--413, 2000.

\bibitem[Ueh14]{Ueha}
Hokuto Uehara.
\newblock Exceptional collections on toric {F}ano threefolds and birational
  geometry.
\newblock {\em Internat. J. Math.}, 25(7):1450072, 32, 2014.

\end{thebibliography}
\bibliographystyle{alpha}
\nocite{M2}
\end{document}